\documentclass[letterpaper, 10pt]{amsart}
\usepackage[all]{xy}        % Commutative Diagrams
\CompileMatrices            % Faster Tex-ing
\UseTips                    % Use Computer Modern Arrowheads
\usepackage[german,english]{babel}   %Foreign language support
\usepackage[bookmarks=true]{hyperref}   % Hyperref in DVI and PDF (like HTML Links between sections)
\usepackage{xspace}
\usepackage{calc}

\theoremstyle{plain}
\newtheorem{theorem}{Theorem}[section]
\newtheorem*{theorem*}{Theorem}
\newtheorem{proposition}[theorem]{Proposition}
\newtheorem{corollary}[theorem]{Corollary}
\newtheorem{lemma}[theorem]{Lemma}

\theoremstyle{definition}
\newtheorem{definition}[theorem]{Definition}

\theoremstyle{remark}
\newtheorem{remark}[theorem]{Remark}
\newtheorem{example}[theorem]{Example}

\newcommand{\enm}[1]{\ensuremath{#1}}          % Shortcuts
\newcommand{\op}[1]{\operatorname{#1}}
\newcommand{\cal}[1]{\mathcal{#1}}

\renewcommand{\bar}[1]{\overline{#1}}

\newcommand{\NN}{\enm{\mathbb{N}}}

\newcommand{\FF}{\enm{\mathbb{F}}}

\newcommand{\Aa}{\enm{\cal{A}}}           % All caligraphy letters easily accessable
\newcommand{\Bb}{\enm{\cal{B}}}

\newcommand{\Ll}{\enm{\cal{L}}}
\newcommand{\Mm}{\enm{\cal{M}}}
\newcommand{\Nn}{\enm{\cal{N}}}

\renewcommand{\phi}{\varphi}        % Dont know how to not loose the original ones???
\renewcommand{\theta}{\vartheta}
\renewcommand{\epsilon}{\varepsilon}

\newcommand{\Spec}{\op{Spec}}

\newcommand{\Hom}{\op{Hom}}

\newcommand{\End}{\op{End}}
\newcommand{\Aut}{\op{Aut}}

\newcommand{\id}{\op{id}}
\newcommand{\dirlim}{\varinjlim}

\newcommand{\tensor}{\otimes}         % Symbols with meaning

\newcommand{\Union}{\bigcup}

\newcommand{\set}[1]{\left\{#1\right\}}
\newcommand{\xn}[1][x]{\enm{#1_1,\ldots,#1_n}}

\renewcommand{\to}[1][]{\xrightarrow{\ #1\ }}

\newcommand{\diff}[1][x]{\textstyle{\frac{\partial}{\partial #1}}}

\newcommand{\usc}[1][m]{\underline{\phantom{#1}}}

\newcommand{\defeq}{\stackrel{\scriptscriptstyle \op{def}}{=}}

\newcommand{\ie}{i.e.,\xspace}           % i.e. in italics and with proper spacing afterwards
\newcommand{\eg}{e.g.,\xspace}           % e.g.   ....
\newcommand{\cf}{{cf.}\ }

\newcommand{\F}[1][]{F^{#1*}}
\newcommand{\T}[1][]{T^{#1*}}
\newcommand{\comment}[1]{}
\newcommand{\eqnref}[1]{(\ref{#1})}

\newenvironment{descrip}[1]
    {\begin{list}{}{
    \settowidth{\labelwidth}{#1}
    \setlength{\leftmargin}{\labelwidth+\labelsep}}}
   {\end{list}}

\begin{document}
\title{The $D$--module structure of $R[F]$--modules}
\author{Manuel Blickle}
\address{Universit\"at Essen, FB6 Mathematik, 45117 Essen, Germany}
\urladdr{www.mabli.org}
\email{manuel.blickle@uni-essen.de}
\thanks{This article was published: Trans. Amer. Math. Soc. 355 (2003), 1647-1668}
\maketitle

\begin{abstract}
Let $R$ be a regular ring, essentially of finite type over a perfect field $k$.
An $R$--module $\Mm$ is called a unit $R[F]$--module if it comes equipped with
an isomorphism $\F[e]\Mm \to \Mm$, where $F$ denotes the Frobenius map on
$\Spec R$, and $\F[e]$ is the associated pullback functor. It is well known
that $\Mm$ then carries a natural $D_R$--module structure. In this paper we
investigate the relation between the unit $R[F]$--structure and the induced
$D_R$--structure on $\Mm$. In particular, it is shown that if $k$ is
algebraically closed and $\Mm$ is a simple finitely generated unit
$R[F]$--module, then it is also simple as a $D_R$--module. An example showing
the necessity of $k$ being algebraically closed is also given.
\end{abstract}

%\tableofcontents

\section{Introduction}
The purpose of this article is to investigate the relationship between
Frobenius actions and differential structure on modules over a regular ring
$R$. Both Frobenius and $D$--module techniques were used with great success in
commutative algebra. For example, in his groundbreaking work, Lyubeznik shows
various finiteness properties of local cohomology modules of a regular ring,
using Frobenius techniques in finite characteristic \cite{Lyub} (see also
\cite{HuSha.LocCohom}) and $D$--modules in characteristic zero
\cite{Lyub.FinChar0}. In fact, these two viewpoints are somewhat reconciled by
noticing that modules with a certain Frobenius action (unit $R[F]$--modules,
see the definition in the next section) carry a natural $D_R$--module
structure. A more careful analysis then shows \cite{Lyub.FinProp,Lyub.InjDim}
that the characteristic zero and $p>0$ proofs are different manifestations of
the same argument.

In this paper we further clarify this connection between the unit
$R[F]$--structure and the induced $D_R$--structure on an
$R$--module $\Mm$. Our main result implies that if $R$ is a
regular ring, essentially of finite type over an
\emph{algebraically closed} field, then a simple, finitely
generated unit $R[F]$--module $\Mm$ is $D_R$--simple. This is a
consequence of the following theorem.
\begin{theorem}\label{thm.main}
    Let $R$ be regular and essentially of finite type over the perfect field $k$.
    Let $\Mm$ be a finitely generated unit $R[F]$--module. Then its
    geometric length as an $R[F]$--module is the same as its geometric
    $D_R$--module length.
\end{theorem}
The term \emph{geometric} refers to the length after tensoring with an
algebraically closed extension field. In addition to correcting an inaccuracy
in the literature (\cite{Lyub}, Remark 5.6a), this result plays an important
part in establishing the finite characteristic analog of the
Kashiwara--Brylinski intersection homology $D_R$--module
\cite{Manuel.PhD,Manuel.int}. For a normal quotient $A=R/I$ of the regular ring
$R$, this module, call it $\Ll(A,R)$, is constructed as the unique simple unit
$R[F]$--submodule of the local cohomology $H^c_I(R)$, using Frobenius
techniques. With the above result one concludes that it is also $D_R$--simple.
Thus $\Ll(A,R)$ is a true analog of Brylinski and Kashiwara's module, which is
characterized as the unique simple $D_R$--submodule of $H^c_I(R)$.

This paper is structured as follows: In Section 2 the basic facts about
$R[F]$--modules are briefly recalled. The main results from \cite{Lyub} are
presented in the form in which they are needed to prove the above statements.
Most proofs are omitted, since they can be found either in \cite{Lyub} or
\cite{Manuel.PhD} with complete detail.

Section 3 continues with background results on $R[F]$--modules, focusing on the
interaction with $D_R$--module theory. In this context a baby version of
Frobenius descent is introduced and applied to derive some elementary
$D_R$--structure results of unit $R[F]$--modules. Then we continue to lay out
all the necessary tools for the proof of Theorem \ref{thm.main}.

Section 4 introduces \emph{geometric $uR[F]$--length} and proves the basic
results about this and related notions of length. Together with the results
from Section 3, we are then able to prove Theorem \ref{thm.main}.

Section 5 contains some examples showing that geometric length possibly
deviates from the length as a unit $R[F]$--module in general. Consequently, the
same examples also give an example of a simple, finitely generated unit
$R[F]$--module which is not $D_R$--simple.

\subsubsection*{Acknowledgements} The results in this article are
part of my dissertation at the Unitversity of Michigan. It is a
great pleasure to thank my advisor Karen Smith for her guidance
and support during my graduate studies.

Thanks go to Matt Emerton for entrusting me with early manuscripts of
\cite{Em.Kis,Em.Kis2}, and to Brian Conrad for numerous comments on my
dissertation, enhancing the content and exposition of this article greatly.

\section{$R[F]$--modules: a brief recall.}\label{sec.recall}
Throughout this paper $R$ always denotes a regular ring. Mostly we will assume
that $R$ is essentially of finite type over a perfect field $k$. The
\emph{(absolute) Frobenius map} on $R$, \ie the ring map sending each element
to its $p$th power, is denoted by $F=F_R$. The associated map on $X=\Spec R$ we
also denote by the same letter $F=F_X$.\footnote{Alternatively we could use the
relative Frobenius, and a similar theory would develop. To keep notation as
simple as possible I chose to stick with the absolute case.}

If $M$ is an $R$--module, then $M^e$ denotes the $R$--$R$--bimodule, which as a
left module is just $M$, but with right structure twisted by the $e$th iterate
of the Frobenius, \ie for $r \in R$ and $m \in M$ one has $m \cdot r =
r^{p^e}m$. With this notation, Psekine and Szpiro's \emph{Frobenius functor} is
defined as $F^{*}(M)=R^1 \tensor M$. Thinking of $F$ as a map on $\Spec R$,
this is just the pullback functor for the Frobenius map. The following
summarizes a few well-known properties of $\F$:
\begin{proposition}\label{prop.basfrob}
\begin{enumerate}
    \item $\F$ is right exact and commutes with direct limits.
    \item If $R$ is regular, then $\F$ is exact and hence commutes with
    finite intersections.
\end{enumerate}
\end{proposition}
Analogously, one defines the higher powers of the Frobenius map and functor.
Obviously, $(F^e)^* = (F^*)^e$, and therefore we denote these higher powers of
the Frobenius functor just by $F^{e*}$. Whenever $e$ is clear from the context,
we denote $p^e$ by $q$.

\begin{definition}
An \emph{$R[F^e]$--module} is an $R$--module $\Mm$ together with
an $R$--linear map
\[
\theta^e: \F[e]\Mm = R^e \tensor \Mm \to \Mm.
\]
If $\theta^e$ is an isomorphism, then $(\Mm,\theta^e)$ is called a
\emph{unit $R[F^e]$--module}.
\end{definition}
By adjointness of extension and restriction of scalars for the Frobenius map on
$\Spec R$, these maps $\theta^e \in \Hom(\F[e]\Mm,\Mm)$ are in one-to-one
correspondence with maps $F^e_\Mm \in \Hom(\Mm,F_*^e\Mm)$. Thus, alternatively,
the $R[F^e]$--module is determined by a map $F^e_\Mm: \Mm \to \Mm$ satisfying
the $p^e$--linearity condition $F^e_\Mm(rm)=r^{p^e}F^e_\Mm(m)$ for all $r \in
R$ and $m \in \Mm$. The relation between the \emph{Frobenius structure}
$\theta^e_{\Mm}: \F[e]\Mm \to \Mm$ and the \emph{Frobenius action} $F^e_{\Mm}:
\Mm \to \Mm$ is illustrated by the commutation of the following diagram:
\begin{equation}\label{dia.frob}
\begin{split}
\xymatrix{
    {R^e \tensor \Mm} \ar@{=}[r] & {\F[e]\Mm} \ar^{\theta^e_\Mm}[d] \\
    {\Mm} \ar^{F^e_R \tensor \id_{\Mm}}[u] \ar^{F^e_\Mm}[r] & {\Mm}
    }
\end{split}
\end{equation}
Concretely, $\theta^e(r \tensor m)=rF^e(m)$ and $F^e(m)=\theta^e(1 \tensor
m)$. For convenience, the subscript ``$\Mm$'' is often omitted on
$\theta^e_\Mm$ and $F^e_\Mm$.

Now it is easy to convince oneself that such a map $F^e_\Mm$ is nothing but an
action of the ring $R[F^e]$ on $\Mm$, where $R[F^e]$ denotes the
(non-commutative) ring which is obtained from $R$ by adjoining the
non-commutative variable $F^e$ to $R$ and forcing the relations $r^{p^e}F^e =
F^{e}r$ of $p^e$--linearity. In other words, an $R[F^e]$--module as defined
above is nothing but a module over the ring $R[F^e]$.

In this sense one defines the category of $R[F^e]$--modules,
$R[F^e]$--mod, as the module category over this ring $R[F^e]$. As
the module category over an associative ring, $R[F^e]$--mod is an
abelian category. It is easily verified that with this definition
a map of $R[F^e]$--modules $(M,\theta^e_M)$ and $(N,\theta^e_N)$
is a map $\phi:M \to N$ of $R$--modules such that
\[
\xymatrix{
    {\F[e]M} \ar^{\F[e]\phi}[r]\ar_{\theta^e_M}[d] & {\F[e]N} \ar^{\theta^e_N}[d] \\
    M\ar^{\phi}[r] & N
    }
\]
commutes. The category of unit $R[F^e]$--modules, $uR[F^e]$--mod, is the full
subcategory whose objects are those $R[F^e]$--modules that are unit. Using that
for a regular ring the Frobenius functor $\F[e]$ is exact, it follows that
$uR[F^e]$--mod is also an abelian category (see \cite[Chapter 2]{Manuel.PhD}
for details).

For all $r > 0$ one has an inclusion of rings, $R[F^{er}] \subseteq R[F^e]$,
yielding the reverse inclusion $R[F^e]\text{--mod} \subseteq
R[F^{er}]\text{--mod}$ of categories. If $(\Mm,\theta^e)$ is an
$R[F^e]$--module, then the Frobenius structure of $\Mm$, viewed as an
$R[F^{er}]$--module, is defined inductively as
$\theta^{er}=\theta^{e(r-1)}\F[e(r-1)](\theta^e)$. It follows that the
categories $\{\, R[F^e]\text{--mod}\,|\, e \in \NN \, \}$ form a directed
system whose limit is denoted by $R[F]$--mod. Thus an $R[F]$--module is an
$R[F^e]$--module for some $e$, not necessarily explicitly
specified.\footnote{In \cite{Manuel.PhD} these categories are denoted by
$R[F^\infty]$--mod and $uR[F^\infty]$--mod respectively. The idea behind the
new notation here is that omission of an exponent on $F$ means an unspecified
exponent and \emph{not} exponent one. Thus $R[F]$--mod seems more adequate to
denote this category, as well as being notationally much more convenient.} It
is easily checked that the category of $R[F]$--modules is also abelian.
Analogously one obtains the category of \emph{unit $R[F]$--modules} as the
directed limit of the categories of unit $R[F^e]$--modules for various $e$.
Again, $uR[F]$--mod is abelian.

\begin{example}
The map $R^e \tensor R \to R$ sending $r' \tensor r \mapsto
r'r^{q}$ makes $R$ into a unit $R[F^e]$--module.

Similarly, a localization $S^{-1}R$ of $R$ is a unit $R[F^e]$--module via
the map $r' \tensor \frac{r}{s} \mapsto \frac{r'r^{q}}{s^{q}}$. The
inverse of this map is given by sending $\frac{r}{s}$ to $rs^{q-1} \tensor
\frac{1}{s}$.
\end{example}

\subsubsection{Functors and base change}
Let $G$ be an additive functor from $R$--mod to $A$--mod, where $A$ is also a
ring of finite characteristic $p$. If $G$ commutes with the Frobenius functor
(\ie we have a natural transformation of functors $\F[e]_A \circ G \cong G
\circ \F[e]_R$), then $G$ naturally extends to a functor from $R[F^e]$--mod to
$A[F^e]$--mod. Indeed, if $(M,\theta^e)$ is a (unit) $R[F^e]$--module, then
\[
    \F[e](G(M)) \cong G(\F[e]((M)) \to[G(\theta^e)] G(M)
\]
defines a (unit) $A[F^e]$--module structure on $G(M)$. The first map is, of
course, the natural isomorphism of functors which was assumed. To verify that
this definition is functorial is straightforward. Since $G$ commuting with
$\F[e]$ implies that $G$ also commutes with the higher powers of $\F[e]$, the
construction just described is indeed a functor on (unit) $R[F]$--modules, \ie
compatible with the inclusion of categories $R[F^e]\text{--mod} \subseteq
R[F^{er}]\text{--mod}$ for $r > 0$.

It follows that for a map of rings $R \to A$, tensoring with $A$ over $R$
is a functor from (unit) $R[F]$--modules to (unit) $A[F]$--modules. All we
need for this is the natural isomorphism of $A$--$R$--bimodules
\[
    A^e \tensor_A A \cong A^e \cong A \tensor_R R^e.
\]
Concretely, if $(\Mm,\theta^e)$ is an $R[F]$--module, then the
$A[F]$--structure on $A \tensor \Mm$ is given by $a' \tensor a \tensor m
\mapsto a'a^q\theta^e(1 \tensor m)$ (equivalently, the Frobenius action is
defined by $F^e_{A \tensor \Mm}=F^e_A \tensor F^e_{\Mm}$).

\subsubsection{Notational conventions and Frobenius.}
If $M$ is an $R$--submodule of the $R[F]$--module $(\Mm, \theta^e, F^e)$, then
$\F[e]M$ is a submodule of $\F[e]\Mm$. Its image in $\Mm$ under $\theta^e$ we
denote by $RF^e(M) = \theta^e(\F[e]M)$, or briefly just by $F^e(M)$. Indeed, by
diagram (\ref{dia.frob}), $RF^e(M)$ is the $R$--submodule of $\Mm$ generated by
the elements $F^e(m)$ for $m \in M$. If $\Mm$ is unit, then $\F[e]M$ is, via
$\theta^e$, isomorphic to $RF^e(M)$. More specifically, on the category of
$R$--submodules of $\Mm$, the two functors $\F[e]$ and $RF^e(\usc)$ are
isomorphic since $RF^e(\usc)=\theta^e \circ \F[e](\usc)$. As it turns out,
working with the more explicit $F^e(\usc)=RF^e(\usc)$ makes many arguments more
transparent, since one never leaves the ambient module $\Mm$. This viewpoint
will be used frequently without further mention.

\subsection{Finitely generated unit $R[F]$--modules}
The strength of $R[F]$--mo\-du\-les lies in bringing some sort of finiteness
to, \textit{a priori}, infinitely generated objects, such as local cohomology
modules. This is simply achieved by enlarging the ring from $R$ to the quite
big, non-commutative ring $R[F^e]$. Exactly the same happens when $D_R$--module
techniques are used; there it is the ring of differential operators which
introduces additional structure.
\begin{definition}
    An $R[F]$--module $(\Mm, \theta^e)$ is called finitely
    generated if it is a finitely generated module over the ring
    $R[F^e]$.
\end{definition}
Note that, if $\Mm$ is generated by some finite subset $S \subseteq \Mm$ as an
$R[F^e]$--module, then, as an $R[F^{er}]$--module, $\Mm$ is generated by the
finite set $S \cup F^e(S) \cup \cdots \cup F^{e(r-1)}(S)$. Thus the finitely
generated $R[F^e]$--modules are a subset of the finitely generated
$R[F^{er}]$--modules. Therefore, the category of finitely generated
$R[F]$--modules is, again, just the limit of the categories of finitely
generated $R[F^e]$--modules for various $e$.

We recall some basic facts from Lyubeznik, but offer a slightly different
viewpoint. The basic construction is that of a generator of a unit
$R[F]$--module.  Let $\phi: M \to \F[e]M$ be an $R$--linear map. Consider the
directed system one obtains by taking higher Frobenius powers of this map. The
limit one obtains
\[
\xymatrix{
    {\phantom{F}\Mm\phantom{F}}\ar@{=}[r] & {\dirlim (M \to[\phi] \F[e]M \to[F^*\phi] \F[2e]M \to \cdots
    \quad)}\ar@{=}[d]\\
    {\F[e]\Mm}\ar@{=}[r]\ar^{\cong}[u] & {\dirlim (\phantom{M \to[\phi] }\F[e]M \to[F^*\phi] \F[2e]M \to \cdots
    \quad)}
}
\]
carries a natural unit $R[F^e]$--module structure as indicated,
taking into consideration that $\F[e]$ commutes with direct
limits.

If a unit $R[F]$--module $(\Mm,\theta^e)$ arises in such a fashion, one calls
$\phi$ a \emph{generator} of $(\Mm,\theta^e)$. If $M$ is finitely generated and
$\phi$ is injective then $M$ is called a \emph{root} of $M$. In this case one
identifies $M$ with its isomorphic image in $\Mm = \dirlim \F[er]M$. Thus a
root of a unit $R[F]$--module $\Mm$ is a finitely generated $R$--submodule $M$
such that $M \subseteq RF^e(M)$ and $\Mm = \bigcup_r RF^{er}(M) = R[F^e]M$.

A key observation is the following proposition \cite{Em.Kis}.
\begin{proposition}
    The unit $R[F]$--module $(\Mm,\theta^e)$ is finitely generated
     if and only if $\Mm$ has a root.
\end{proposition}
\begin{proof}
The ``only if'' direction is easy, since a finite set of $R$--module generators
of a root $M$ of $\Mm$ generates $\Mm$ as an $R[F^e]$--module.

Conversely, let $M'$ be the $R$--module generated by some finitely
many $R[F^e]$--module generators of $\Mm$. In other words,
\begin{equation}\label{l.1}
  R[F^e]M'=\sum_{n=0}^{\infty} F^{ne}(M') = \Mm.
\end{equation}
Since $\theta^e$ is an isomorphism, $\Mm = \theta(R^e \tensor \Mm) = F^e(\Mm)$.
Applying $F^e$ to (\ref{l.1}), we get
\[
    \Mm = F^e(\Mm) = F^e(\sum_{n=0}^\infty F^{ne}(M')) = \sum_{n=1}^\infty F^{ne}(M').
\]
Since $M'$ was finitely generated, it is contained in a finite part of the
above sum, say $M' \subseteq \sum_{n=1}^m F^{ne}(M')$. Now we set $M =
\sum_{n=0}^{m-1} F^{ne}(M')$, and we see right away that $M \subseteq F^e(M)$.
Repeated application of $F^e(\usc)$ yields the following sequence of
inclusions:
\[
    M \subseteq F^e(M) \subseteq F^{2e}(M) \subseteq F^{3e}(M)
    \subseteq \cdots,
\]
whose union is $\Mm$ since $M$ contains $M'$. Thus $M$ is the
desired root of $\Mm$.
\end{proof}
This proposition allows for easy proofs of the following theorem
of \cite{Lyub}, Theorem 2.8:
\begin{theorem}
    The category of finitely generated unit $R[F]$--modules is an
    abelian subcategory of the category of $R[F]$--modules that is
    closed under extensions.

    All finitely generated unit $R[F]$--modules have the ascending
    chain condition in the category of unit $R[F]$--modules.
\end{theorem}
\begin{proof}
Everything besides the fact that the kernels are finitely generated easily
follows using either the exactness of $\F[e]$ or basic facts about categories
of finitely generated modules over an associative ring (they are abelian and
closed under extensions). Thus it remains to show that a unit
$R[F^e]$--submodule $\Nn$ of a finitely generated unit $R[F^e]$--module $\Mm$
is also finitely generated. To see this, let $M$ be a root of $\Mm$, \ie $M
\subseteq RF^e(M)$ and $M=R[F^e]M$. Then
\[
    M \cap \Nn \subseteq RF^e(M) \cap \Nn = RF^e(M) \cap RF^e(\Nn) = RF^e(M \cap \Nn)
\]
and $R[F^e](M \cap \Nn)=R[F^e]M \cap R[F^e]\Nn = \Mm \cap \Nn = \Nn$. Thus $M
\cap \Nn$ is a root of $\Nn$, and therefore $\Nn$ is finitely generated as an
$R[F^e]$--module. In this argument the fact (Proposition \ref{prop.basfrob})
that $RF(\usc)$ commutes with finite intersection was used repeatedly.

To see that a finitely generated unit $R[F]$--module has ACC, note that we just
showed that a unit $R[F]$-submodule $\Nn \subseteq \Mm$ is determined by its
root $N$, the intersection of $\Nn$ with a fixed root $M$ of $\Mm$. Thus
increasing chains of unit $R[F]$--submodules of $\Mm$ correspond, by
intersecting with the root $M$, to certain chains of submodules of $M$. The
latter stabilize because $M$ is finitely generated as an $R$--module, and thus
so do the former.
\end{proof}
Already the fact that $uR[F]$--mod has ACC is very useful for proving
finiteness statements for local cohomology modules, such as the finiteness of
the set of associated primes or the finiteness of the Bass numbers. But even
more is true. The following is Theorem 3.2 of \cite{Lyub}; we recall it here
without proof.
\begin{theorem}\label{thm.RFfinite}
    Let $R$ be a regular finitely generated algebra over a regular local ring.
    Then, in the category of unit $R[F^e]$--modules, the finitely generated ones
    have finite length.
\end{theorem}
It is an open problem whether this theorem is true without the
assumptions on $R$.

\section{Frobenius action and differential structure}

The connection between unit $R[F]$--modules and differential
operators originates in the following description of the ring of
differential operators over a ring of finite characteristic.

\begin{proposition}\label{prop.DasEnd}
    Let $R$ be a finitely generated algebra over its subring $R^p$ of $p$th powers. Then
    \[
        D_R= \Union \End_{R^{p^e}}(R)
    \]
    is the ring of $k$--linear differential operators on $R$.
\end{proposition}
For our purpose this proposition can be taken as the definition of the
ring of differential operators, therefore eliminating the need of
reviewing its proof, which can be found in \cite{Yeku}. The sets
$\End_{R^{p^e}}(R)$ are called the differential operators of level $e$ and
are denoted by $D^{(e)}_R$. The interested reader will find a gentle
introduction to differential operators in finite characteristic in
\cite{Manuel.PhD}; for the ultimate account see
\cite{Ber.OpDiff,Ber.FrobDesc}.

With this description of the ring of differential operators the action of $D_R$
on a unit $R[F]$--module becomes fairly straightforward. First note that
$\End_{R^{p^e}}(R)$ can be identified with $\End_{\text{mod--}R}(R^e)$, the
right $R$--module endomorphisms of $R^e$. For big enough $r$, a given
differential operator $\delta$ lies in the set $D^{(er)}_R =
\End_{\text{mod--}R}(R^{er})$. For such $r$ the action of $\delta$ on a unit
$R[F^e]$--module $(M,\theta^e)$ is given by the dashed arrow of the following
diagram:
\[
\xymatrix@C=3pc{ M \ar_{(\theta^{er})^{-1}}[d]\ar^{\delta\cdot}@{-->}[r] & M  \\
           {R^{er} \tensor M} \ar^{\delta \tensor \id_M}[r] & {R^{er} \tensor M} \ar_{\theta^{er}}[u]
            }
\]
To check that this is independent of the chosen $r$ is
straightforward. This also shows that one gets the same
$D_R$--structure whether one considers $\Mm$ as a unit
$R[F^e]$--module or a unit $R[F^{er}]$--module, for some $r>0$.
This implies that the process of equipping a unit $R[F^e]$--module
with an underlying $D_R$--structure is compatible with the
inclusion of categories $uR[F^e]\text{--mod} \subseteq
uR[F^{er}]$, thus is well defined as a functor from $uR[F]$--mod
to $D_R$--mod. The next lemma summarizes some properties of this
$D_R$--structure.
\begin{lemma}\label{lem.DRFProp}
\begin{enumerate}
\item
    The $D_R$--structure on $R$ induced by the canonical unit
    $R[F]$--structure is the canonical $D_R$--structure on $R$.
\item
    The described process of equipping a unit $R[F]$--module with a
    $D_R$--module structure is an exact functor which commutes with
    localization.
\end{enumerate}
\end{lemma}
For more details and proves of these (and the following) statements, the reader
should consult \cite{Manuel.PhD}, Chapter 3, or \cite{Lyub}, Section 5.

\subsection{Frobenius descent}
The interplay between Frobenius structures and $D_R$--modules becomes apparent
by the simple, but powerful, fact that the Frobenius functor is an equivalence
of the category of $D_R$--modules with itself. Ultimately, this is a
consequence of the description of $D_R$ as in Proposition \ref{prop.DasEnd}
together with the observation that all the endomorphism rings
$D^{(e)}_R=\End_{\text{mod--}R}(R^e)$ are Morita equivalent to $R$ itself.
This, on the other hand, is easily seen, since, for a regular $R$, $R^e$ is
locally a free right $R$--module and therefore $\End_{\text{mod--}R}(R^e)$ is
just a matrix ring over $R$, thus Morita equivalent to $R$. Then one checks
that this Morita equivalence between $R$ and $\End_{\text{mod--}R}(R^e)$ is, in
fact, given by the Frobenius functor $\F[e] \cong R^e\tensor \usc$. One obtains
the following theorem.
\begin{theorem}[Frobenius descent]
    Let $R$ be regular and $F$--finite. The Frobenius functor $\F[e]$ is
    an equivalence of categories between the category of
    $R$--modules and $D^{(e)}_R$--modules. The functor inverse to $\F[e]$
    is given by
    \[
        T^{e*} \defeq \Hom_{\text{mod--}R}(R,R^e) \tensor_{D^{(e)}_R} \usc
    \]
    The functors $\F[e]$ and $T^{e*}$ induce an auto-equivalence of the category
    of $D_R$--modules.
\end{theorem}
This result appears, much generalized, in Berthelot \cite{Ber.FrobDesc} under
the name of Frobenius descent. Similar versions are used by S.P. Smith
\cite{SmithSP.diffop,SmithSP:DonLine}, B. Haastert
\cite{Haastert.DiffOp,Haastert.DirIm} and R. B{\o}gvad \cite{Bog.DmodBorel}; also
Lyubeznik implicitly uses Frobenius descent to prove the following theorem.
\begin{theorem}[\cite{Lyub}, Theorem 5.6]\label{thm.Drfinite}
    Let $R$ be regular and $F$--finite. Then, a finitely generated unit
    $R[F^e]$--module that has finite length as a unit $R[F^e]$--module also
    has finite length as a $D_R$--module.

    In particular, if $R$ is a finitely generated algebra
    over a regular local ring, then all finitely generate unit
    $R[F^e]$--modules have finite $D_R$--module length.
\end{theorem}
For the proof, please refer to \cite{Lyub}. In the next section
part of the argument will be given to see explicitly how Frobenius
descent enters.

If $\Mm$ is a $D_R$--module, the $D_R$--structure that $\F[e]\Mm$ obtains via
Frobenius descent is described as follows. Since $\Mm$ is a $D_R^{(e)}$--module
for all $e$, we can write, by Frobenius descent $\Mm \cong \F[e]\T[e]\Mm$,
where $\T[e]\Mm$ is an $R$--module. Thus $\F\Mm \cong \F[(e+1)]\T[e]\Mm$
obtains a $D_R^{(e+1)}$--module structure by having $\delta$ act via $\delta
\tensor \id_{\T[e]\Mm}$. One can check (similarly as one checks the
well-definedness of the induced $D_R$--action from a unit $R[F]$--module) that
the $D_R^{(e)}$--structures on $\F\Mm$ that one obtains in this way for all $e$
are compatible and define a $D_R$--structure on $\Mm$.

\subsubsection{Frobenius descent for unit $R[F]$--modules}
As we have noted before, for an $R$--submodule $N$ of a unit $R[F^e]$--module
$(\Mm,\theta^e,F^e)$, the image $\F[e]N$ of $N$ under the Frobenius functor
finds a concrete description as the $R$--submodule $F^e(N)$ of $\Mm$.
Similarly, for a $D_R$--submodule $\Nn$ of $\Mm$ the $D_R$--module $\T[e]\Nn$
can be realized as a $D_R$--submodule $T^e(\Nn)$ of $\Mm$ as follows: On the
category of $D_R$--submodules of $\Mm$, the functor $\T[e]_\Mm\usc$ can be
identified with $T^e_\Mm(\usc) \defeq \T[e]_\Mm(\theta^{-1}(\usc))$. Together
with the natural identification $\T[e]\F[e]\Mm \cong \Mm$ this makes $T^e(\Nn)$
a submodule of $\Mm$. Clearly, $F^e \circ T^e = T^e \circ F^e$ are the identity
functor on $D_R$--submodules of $\Mm$. A more careful investigation shows that
$T^e_\Mm(\Nn)=(F_{\Mm}^e)^{-1}(\Nn)$ as submodules of $\Mm$. Since we don't
need this description in what follows, its proof is omitted and the reader is
referred to \cite{Manuel.PhD}, Chapter 3, for details.

When working with $D_R$--submodules of a fixed unit $R[F]$--module, it is
advantageous to use the un-starred variants of the functors $\T[e]$ and
$\F[e]$, since this makes the arguments much more transparent.

\subsection{$D_R[F]$--modules} If $(\Mm,\theta^e)$ is a unit
$R[F^e]$--module, then  $\F[e]\Mm$ carries a natural $D_R$--structure by
Frobenius descent, as well as by being a unit $R[F^e]$--module
$(\F[e]\Mm,\F[e](\theta^e))$. To show that these two $D_R$--structures are the
same comes down to checking that $\theta^e$ is $D_R$--linear with respect to
the $D_R$--structure on $\F[e]\Mm$ coming from Frobenius descent. This
motivates us to define a $D_R[F^e]$--module in analogy with $R[F^e]$--modules
as a $D_R$--module $\Mm$ together with a $D_R$--linear map $\theta^e: \F[e]\Mm
\to \Mm$. Of course, the $D_R$--structure on $\F[e]\Mm$ is the one coming from
Frobenius descent.

\begin{definition}
    A \emph{$D_R[F^e]$--module} is a
    $D_R$-module $\Mm$ together with a $D_R$-linear map
    \[
        \theta^e_\Mm: \F[e]_R \Mm \to{} \Mm.
    \]
    In other words, a $D_R[F^e]$--module is an $R[F^e]$--module
    $(\Mm,\theta^e)$ that carries a $D_R$--structure such that $\theta^e$ is
    $D_R$--linear. $(\Mm,\theta^e)$ is called \emph{unit} if $\theta^e$
    is an isomorphism.
\end{definition}
As just elaborated, unit $R[F^e]$--modules are $D_R[F^e]$--modules. Conversely,
by forgetting the $D_R$--structure, every unit $D_R[F^e]$--module is a unit
$R[F^e]$--module. Thus, for the unit case this does not lead to anything new.
One reason for working with the more complicated category of
$D_R[F^e]$--modules is the following proposition.
\begin{proposition}\label{prop.DRfulluRF}
    Let $R$ be regular and $F$--finite. Then a $D_R[F^e]$--submodule of
    a finitely generated unit $R[F^e]$--module is also unit (and finitely
    generated).
\end{proposition}
Note that the equivalent statement for $R[F^e]$--modules is not true. The
example of an ideal $I$ of $R$, which is an $R[F^e]$--submodule but not unit
shows this nicely.
\begin{proof}[Proof of \ref{prop.DRfulluRF}]
We have to show that a $D_R[F^e]$--submodule of a finitely generated unit
$R[F^e]$--module is also unit. If $\Nn$ is such an $F^e$--stable
$D_R$--submodule of a unit $R[F^e]$--module $\Mm$, then $F^e(\Nn) \subseteq
\Nn$. Applying $T^e$ and using its defining property as the inverse functor of
$F^e$, we see that $\Nn \subseteq T^e(\Nn)$. Iterating, we get an increasing
chain of $R$--modules
\[
        \Nn \subseteq T^e(\Nn) \subseteq T^{2e}(\Nn) \subseteq \ldots.
\]
Intersecting this chain with a root $M$ of $\Mm$ yields a chain of submodules
of $M$ which, since $M$ is a finitely generated $R$--module, must stabilize.
Let $N \defeq T^{er}(\Nn) \cap M=T^{e(r+1)}(\Nn)\cap M= \ldots$ be the stable
member. Using $F^{er}(T^{er}(\Nn))=\Nn$, we see that
$F^e(N)=F^e(T^{e(r+1)}(\Nn) \cap M)=T^{er}(\Nn)\cap F^e(M) \supseteq
T^{er}(\Nn) \cap M = N$. Applying $F^{er}$ for all $r$ to this inclusion, we
get another increasing sequence:
\[
    N \subseteq F^e(N) \subseteq F^{2e}(N) \subseteq \ldots
\]
Let $\Ll$ be its limit. Since $\Ll$ arises as the increasing union of the
Frobenius powers of a single submodule, it is obviously a unit submodule of
$\Nn$, \ie $F^e(\Ll)=\Ll \subseteq \Nn$. For the converse inclusion let $n \in
\Nn$. For all sufficiently large $r \geq 0$ we have $n \in F^{er}(M)$. For such
$r$ also $N=T^{er}(\Nn) \cap M$, and thus $F^{er}(N)=\Nn \cap F^{er}(M)$. Thus
$n \in F^{er}(N)$, and therefore $n \in \Ll$, since $\Ll$ is the increasing
union of all $F^{er}(N)$.
\end{proof}
\begin{remark}
    Analogously to the case of $R[F^e]$--modules, one can show that
    $D_R[F^e]$--modules are just modules over an appropriate ring
    $D_R[F^e]$. In the case that the ring $R[F^e]$ is in fact a subring of
    $\End_k(R)$ we can think of this ring
    $D_R[F^e]$ as the subring of $\End_k(R)$ generated by $R[F^e]$ and
    $D_R$. In general one can define the ring $D_R[F^e]$ to be $R[F^e]
    \tensor_R D_R$, and then equip this tensor product with an appropriate
    ring structure. This is done in \cite{Em.Kis2}, where many other
    interesting properties in this context are shown. For example, they
    show the following (at first) surprising result:
    \begin{proposition}
        Let $R$ be a regular ring, essentially of finite type over a perfect field
        $k$. Let $M$ be a unit $R[F^e]$--module. The following are
        equivalent.
        \begin{enumerate}
        \item $M$ is finitely generated as a $D_R[F^e]$--module.
        \item $M$ is finitely generated as an $R[F^e]$--module.
        \item $M$ is finitely generated as a $D_R$--module.
        \end{enumerate}
    \end{proposition}
    The proof of this is a clever application of Frobenius descent
    together with \cite{Lyub}, Theorem 5.6.
\end{remark}

\subsection{$D_R$--semisimple unit
$R[F]$--modules.}\label{sec.semisimp}

An object $\Mm$ in an abelian category is called semisimple if every subobject
is a direct summand. If $\Mm$ also has finite length (in the category), this is
equivalent to every simple subobject being a direct summand. In turn, this can
be shown to be the same as $\Mm$ decomposing into a finite direct sum of simple
objects. The first result indicating the importance of the notion of
semisimplicity in the context of unit $R[F]$--modules and $D_R$--modules is the
following proposition which is implicit in the proof of \cite{Lyub}, Theorem
5.6.
\begin{proposition}\label{prop.RFsemisimpD}
    Let $R$ be regular and $F$--finite.
    Let $(\Mm,\theta^e)$ be a simple finitely generated unit $R[F^e]$--module.
    Then $\Mm$ is semisimple as a $D_R$--module.
\end{proposition}
\begin{proof}
We show that every simple $D_R$--submodule $\Nn$ of $\Mm$ is a direct
$D_R$--module summand. By Frobenius descent, $F^e\Nn$ is a simple
$D_R$--submodule of $\Mm$. Repeatedly applying $F^e(\usc)$, we get a series
$\{\, F^{re}\Nn\, \}$ of simple $D_R$--submodules of $M$. Let $r$ be the first
time such that the intersection
\[
    (\Nn + F^e\Nn + F^{2e}\Nn + \ldots + F^{re}\Nn) \cap
    F^{(r+1)e}\Nn
\]
is nonempty. Since $\Mm$ has finite $D_R$--module length by Theorem
\ref{thm.Drfinite}, such $r$ exists. For this $r$, the sum on the left is
direct; in particular, $\Nn$ is a direct summand of $\Mm'
\defeq \Nn \oplus F^e\Nn \oplus F^{2e}\Nn \oplus \ldots \oplus
F^{re}\Nn$. Furthermore, $F^{(r+1)e}\Nn$ is contained in $\Mm'$, by simplicity
of $F^{(r+1)e}\Nn$. Then $F^e\Mm' \subseteq \Mm'$, and therefore $\Mm'$ is a
$D_R$-- and $R[F^e]$--submodule of $\Mm'$. By Proposition \ref{prop.DRfulluRF},
$\Mm'$ is a \emph{unit} $R[F^e]$--submodule of $\Mm$. By simplicity of $\Mm$ as
a unit $R[F^e]$--module it follows that $\Mm'=\Mm$. Thus $\Nn$ is a direct
summand of $\Mm$.
\end{proof}
\begin{proposition}
    Let $R$ be regular and $F$--finite, and let
    $(\Mm,\theta^e)$ be a finitely generated unit $R[F^e]$--module that is
    semisimple as a $D_R$--module. For some $r > 0$, all
    $D_R$--isotypic components of $M$ are $R[F^{er}]$--submodules of
    $\Mm$.
\end{proposition}
\begin{proof}
Let $\Mm = \Mm_1 \oplus \ldots \oplus \Mm_s$ be the decomposition
of $\Mm$ into its isotypic components as a $D_R$--module, \ie
$\Mm_i \cong \Nn_i^{\oplus n_i}$ for distinct, simple
$D_R$--modules $\Nn_i$.

By Frobenius descent, $\F[e]\Nn_i$ are again distinct (pairwise non-isomorphic)
simple $D_R$--modules. Since $\F[e]$ commutes with finite direct sums, it
follows that $\F[e]\Mm_1 \oplus \ldots \oplus \F[e]\Mm_s$ is the isotypic
decomposition of $\F[e]\Mm$. Since $\theta^e$ is a $D_R$--module isomorphism,
it maps an isotypic component $\F[e]\Mm_i$ isomorphically onto another isotypic
component $\Mm_{\sigma(i)}$ for some permutation $\sigma$ of the index set
$\{\,1, \ldots,s \, \}$. For $r$, the order of $\sigma$ (\ie $\sigma^r = \id$),
it follows that $\theta^{er}$ is an isomorphism of $\F[er]\Mm_i$ onto $\Mm_i$
for all $i$. Thus, the $D_R$--isotypic decomposition of $\Mm$ is, in fact, a
decomposition of $\Mm$ as a unit $R[F^{er}]$--module.
\end{proof}
\begin{corollary}\label{cor.IsotypicAreSub}
    Let $R$ be regular and $F$--finite, and let $\Mm$ be a finitely
    generated unit $R[F]$--module that is
    semisimple as a $D_R$--module. Each $D_R$--isotypic component of
    $\Mm$ is an $R[F]$--submodule of $\Mm$.

    In particular, a simple, finitely generated unit $R[F]$--module
    $\Mm$ is $D_R$--isotypic, \ie as a $D_R$--module, $\Mm \cong
    \Nn^{\oplus n}$ for some simple $D_R$--module $\Nn$.
\end{corollary}
As a consequence, we get that a simple $D_R$--submodule $\Nn$ of a unit
$R[F]$--module $\Mm$ carries also a unit $R[F]$--structure, although the
inclusion $\Nn \subseteq \Mm$ is only $R[F]$--linear if $F^e_{\Mm}(\Nn)
\subseteq \Nn$, which is in general not the case.
\begin{proposition}\label{prop.SimpleSub}
    Let $R$ be regular and $F$--finite. Let $\Mm$ be a finitely generated
    unit $R[F]$--module that has finite unit $R[F]$--module length.
    Let $\Nn$ be a $D_R$--sub-quotient of $\Mm$.
    Then $\Nn$ carries a unit $R[F]$--module structure.
\end{proposition}
\begin{proof}
Since the category of unit $R[F]$--modules is closed under extensions, we can
reduce, by induction on the $D_R$--module length of $\Nn$ (which is finite by
Theorem \ref{thm.Drfinite}), to the case of $\Nn$ being $D_R$--simple.

Using that $\Mm$ has finite length as a unit $R[F]$--module and as a
$D_R$--module, one reduces to the case that $\Mm$ itself is simple as a unit
$R[F]$--module. Then, by Corollary \ref{cor.IsotypicAreSub}, $\Mm$ is
semisimple as a $D_R$--module and $\Nn$--isotypic, \ie as a $D_R$--module, $\Mm
\cong \Nn^{\oplus n}$. Let $\theta^e$ be the structural morphism of $\Mm$. Let
$\pi_i: \Mm \to \Nn$ be the projection onto the $i$th direct summand. The
composition
\[
    \F[e]\Nn \subseteq \F[e]\Mm \to[\theta^e] \Mm \to[\pi_i] \Nn
\]
is nonzero for at least one index $i$, since $\theta^e$ is an isomorphism. As a
$D_R$--linear map between simple modules, it must be an isomorphism. Thus
$\F[e]\Nn \cong \Nn$ and $\Nn$, is a unit $R[F]$--module.
\end{proof}

Note that the constructed unit $R[F]$--structure on $\Nn$ is by no means
unique;\footnote{In general, the unit $D_R[F^e]$--structures of a $D_R$--simple
module $\Nn$ correspond to $\End_{D_R}(\Nn)$, which is expected to be a rather
controllable set if $\Nn$ is $D_R$--simple, \cf \ref{sec.assump}.} neither are
the maps which express $\Nn$ as a sub-quotient of $\Mm$ maps of
$R[F]$--modules.

\subsection{Connection to vector-spaces with Frobenius
action}\label{sec.VS}

By definition, a $D_R[F^e]$--module structure on $\Mm$ is a $D_R$--linear map
$\theta^e: \F[e]\Mm \to \Mm$ where $\F[e]\Mm$ carries its natural
$D_R$--structure given by Frobenius descent. Thus, the set of
$D_R[F^e]$--module structures on $\Mm$ can be identified with
$\Hom_{D_R}(\F[e]\Mm,\Mm)$. If $\Mm$ admits a unit $D_R[F^e]$--structure
$\theta^e:\F[e]\Mm \to[\cong] \Mm$, this induces an isomorphism
$\Hom_{D_R}(\F[e]\Mm,\Mm)\cong \End_{D_R}(\Mm)$.
\begin{lemma}\label{lem.DFstructs}
    Let $(\Mm,\theta^e)$ be a unit $D_R[F^e]$--module. Every
    $D_R[F^e]$--structu\-re ${\theta'}^e$ on $\Mm$ can be written as
    ${\theta'}^e=\phi \circ \theta^e$ for some unique $\phi \in
    \End_{D_R}(\Mm)$.
\end{lemma}
\begin{proof}
The map ${\theta'}^e \mapsto {\theta'}^e \circ (\theta^e)^{-1}$ is inverse
to the map $\phi \mapsto \phi \circ \theta^e$.
\end{proof}
Furthermore, ${\theta'}^e$ defines a unit $D_R[F]$--structure if and only if
the associated $\phi = {\theta'}^e \circ (\theta^e)^{-1}$ in $\End_{D_R}(\Mm)$
is a $D_R$--module automorphism of $\Mm$. Thus the unit $D_R[F^e]$--structures
are in one-to-one correspondence with $\Aut_{D_R}(\Mm)$.

Also note that the Frobenius actions corresponding to $\theta^e$ and
${\theta'}^e$ on $\Mm$ are related by ${F'}^e = \phi \circ F^e$.

Now assume that $\Mm$ is $D_R$--isotypic. Then $\Mm \cong V \tensor_k \Nn$ for
some finite dimensional $k$--vector-space $V$ and a simple $D_R$--module $\Nn$.
In this way of writing things a differential operator $\delta$ acts as
$\delta(v \tensor n)=v \tensor \delta(n)$. This is well defined since $k$ is
perfect.

In order to study the Frobenius actions on $\Mm=V \tensor_k \Nn$ via the
Frobenius actions on $V$, we assume that $\End_{D_R}(\Nn) = k$. This allows us
to identify the $D_R$--module endomorphisms of $\Mm$ with the $k$--linear
endomorphisms of $V$.

\subsubsection{Discussion of the assumption
$\End_{D_R}(\Nn)=k$}\label{sec.assump} In characteristic zero, under the
assumption that $R$ is a regular ring, essentially of finite type over $k$,
Quillen's lemma \cite{Quill} shows that the $D_R$--module endomorphism set of a
simple $D_R$--module is algebraic over $k$. Thus, if $k$ is algebraically
closed, $\End_{D_R}(\Nn)=k$. Quillen's proof exploits the fact that $D_R$, in
characteristic zero, is a finitely generated algebra over $k$. This fails in
finite characteristic. Nevertheless, one has the following result due to
Dixmier \cite{Dix} for convenience we also recall the proof.
\begin{lemma}\label{lem.EndNisAlgebraic}
    Let\/ $R$ be a $k$--algebra such that the cardinality of $k$ is strictly
    bigger than the cardinality of a $k$--basis of $R$.
    If\/ $\Nn$ is a simple\/ $D_R$--module, then\/ $\End_{D_R}(\Nn)$ is
    algebraic over\/ $k$.
\end{lemma}
\begin{proof}
Let $\kappa$ be the cardinality of the $k$--basis of $R$ and $\kappa'$ the
strictly bigger cardinality of $k$. Since $\Nn$ is simple, we have $D_Rn=\Nn$
for some (every nonzero) $n \in \Nn$. Thus every $\phi \in \End_{D_R}(\Nn)$ is
determined by its value on $n$. Since $D_R$ is at most $\kappa$--dimensional
over $k$, so is $\Nn=D_Rn$, and thus $\End_{D_R}(\Nn)$ is also at most
$\kappa$--dimensional over $k$. Therefore, for any fixed $\phi \in
\End_{D_R}(\Nn)$ (say $\phi \not\in k$), the set $\set{(\phi + \lambda)^{-1}\
|\ \lambda \in k\ }$ has cardinality $\kappa'>\kappa$, thus must be linearly
dependent (we use that $\End_{D_R}(\Nn)$ is a division ring by Shur's lemma). A
relation of linear dependence among some finitely many $(\phi +
\lambda_i)^{-1}$ gives, after clearing denominators, an algebraic relation for
$\phi$. Clearing denominators works just as in the commutative case, since all
$(\phi + \lambda_i)^{-1}$ commute with each other.
\end{proof}
It is an interesting open problem whether an analog of Quillen's lemma holds in
finite characteristic.

For the rest of this section we assume that $\End_{D_R}(N)=k$. In the next
section we will be able to put ourselves in a situation of Lemma
\ref{lem.EndNisAlgebraic} so that we are able to apply the following results.
\begin{lemma}\label{lem.EndDEndk}
    Let $\Mm \cong V \tensor_k \Nn$ for a simple $D_R$--module
    $\Nn$ such that $\End_{D_R}(\Nn)=k$. Then $\End_{D_R}(\Mm) \cong
    \End_k(V)$.

    Tensoring with $\Nn \tensor \usc$ gives a one-to-one
    correspondence between the $k$--vector subspaces of $V$ and the
    $D_R$--submodules of $\Nn \tensor V \cong \Mm$.
\end{lemma}
\begin{proof}
After the choice of a basis for $V$, the ring $\End_{D_R}(\Mm)$ is the matrix
algebra over $\End_{D_R}(\Nn)=k$ of size $\dim_k V$. This identifies
$\End_{D_R}(\Mm)$ with $\End_k(V)$. Given a $\phi_k \in \End_k(V)$, the
corresponding map in $\End_{D_R}(\Mm)$ is $\id_\Nn \tensor \phi_k$.

Let $\Mm'$ be a $D_R$--submodule of $\Mm$. Since $\Mm$ is semisimple we find a
$D_R$--submodule $\Mm''$ such that $\Mm' \oplus \Mm'' \cong \Mm$. Then $\Mm'$
is the kernel of the endomorphism $\phi: \Mm \to[\pi]\Mm'' \subseteq \Mm$,
where $\pi$ is the projection onto the direct summand $\Mm''$. Thus, by the
first part, $\phi = \id_\Nn \tensor \phi_k$ for some $\phi_k$ in $\End_k(V)$.
Then clearly $\Mm' = V' \tensor \Nn$, with $V' = \ker \phi_k$.
\end{proof}
I want to extend the last lemma so that it incorporates Frobenius
operations. Then it is possible to reduce questions about unit
$R[F]$--submodules of $M$ (\eg are there any nontrivial ones?) to the
equivalent questions about $k[F]$--subspaces of $V$. The main observation
is:
\begin{proposition}\label{prop.MreducetoV}
    Let\/ $\Mm=V \tensor_k \Nn$ be as in the last lemma. Assume that $\Mm$
    is a unit\/ $D_R[F^e]$--module with Frobenius action\/ $F^e_{\Mm}$. Then there
    is a Frobenius action\/ $F^e_V$ on\/ $V$ such that the\/
    $D_R[F^e]$--submodules of\/ $\Mm$ are in one-to-one correspondence with the\/
    $k[F^e]$--submodules of\/ $V$.
\end{proposition}
\begin{proof} First observe that by Proposition
\ref{prop.SimpleSub}, $\Nn$ carries some (not necessarily unique)
 unit $D_R[F^e]$--module structure. Denote the
corresponding Frobenius action by $F^e_\Nn$. The choice of a basis of $V$
equips $V$ with a unit $k[F^e]$--structure by letting $F^e$ act as the identity
on the basis and extending $p^e$-linearly. This Frobenius action we denote by
${F'}^e_V$. Then ${F'}^e_\Mm
\defeq F^e_\Nn \tensor {F'}^e_V$ defines a unit $D_R[F^e]$--structure.
The corresponding Frobenius structure on $\Mm$ is given by
\[
    (R^e \tensor_R \Nn) \tensor_k V \to \Nn^e \tensor_k V \cong \Nn \tensor_k
(k^e \tensor_k V) \to \Nn \tensor_k V,
\]
where the first map is the unit $R[F^e]$--structure on $\Nn$ and
the last is the unit $k[F^e]$--structure on $V$. Since these both
are isomorphisms, so is the composition.

By Lemma \ref{lem.DFstructs}, we can express $F^e_\Mm=\phi \circ {F'}^e_\Mm$
for some $\phi \in \Aut_{D_R}(\Mm)$. By Lemma \ref{lem.EndDEndk} we can write
$\phi = \id_\Nn \tensor \phi_V$ for some $k$--vector-space automorphism
$\phi_V$ of $V$. Denoting the corresponding unit $k[F^e]$--structure on $V$ by
$F^e_V \defeq \phi_V \circ {F'}^e_V$, one easily verifies that $F^e_\Mm =
F^e_\Nn \tensor F^e_V$. Indeed,
\[
    F^e_\Nn \tensor F^e_V = F^e_\Nn \tensor (\phi_V \circ {F'}^e_V) = \phi \circ
    (F^e_\Nn \tensor {F'}^e_V) = \phi \circ {F'}^e_\Mm = F^e_\Mm.
\]
With $F^e_V$ we have constructed the desired Frobenius action on $V$. Since
$F^e_\Mm=F^e_\Nn \tensor F^e_V$, we have $F^e_\Mm(\Nn \tensor V')=\Nn \tensor
F^e_V(V')$. Therefore, a $D_R$--submodule of the form $\Nn \tensor V'$ of $\Mm$
is stable under $F^e_\Mm$ if and only if the subspace $V'$ of $V$ is stable
under $F^e_V$. By Lemma \ref{lem.EndDEndk} every $D_R$--submodule of $\Mm$ is
of this form.
\end{proof}
Note that by Proposition \ref{prop.DRfulluRF}, the $D_R[F]$--submodules of
$\Mm$ are exactly the unit $R[F]$--submodules. One gets the following
corollary.
\begin{corollary}
With the notation as in the last proposition, the unit
$R[F^e]$--submodules of $(\Mm,F^e_{\Mm})$ are in one-to-one correspondence
with the $k[F^e]$--submodules of $(V,F^e_V)$.

In particular, $\Mm$ is a simple unit $R[F]$--module if and only
if $V$ is a simple $k[F]$--module.
\end{corollary}

\section{Length of unit $R[F]$--modules}
The length of an object $\Mm$ in an abelian category (the length of the longest
chain of proper inclusions of objects in the category starting with zero and
ending with $\Mm$) is denoted by small $l_*$ decorated by a modifier pointing
out the category. For example, if $\Mm$ is a finitely generated unit
$R[F^e]$--module, we denote
\begin{eqnarray*}
    l_{uR[F^e]}(\Mm) &=& \text{ the length of the unit $R[F^e]$--module
    $\Mm$,} \\
    l_{uR[F^{er}]}(\Mm) &=& \text{ the length of $\Mm$ as a unit
    $R[F^{er}]$--module} \\
    l_{D_R}(\Mm) &=& \text{ the length of $\Mm$ as a $D_R$--module,} \\
    l_{uR[F]}(\Mm) &=& \text{ the length of $\Mm$ as a unit $R[F]$--module.}
\end{eqnarray*}
One easily concludes the following proposition.

\begin{proposition}\label{prop.lenghtsFinite} Let $R$ be essentially of
finite type over the perfect field $k$. Let $\Mm$ be a finitely generated
unit $R[F^e]$--module. Then there is an $r
> 0$ such that
\[
l_{uR[F^e]}(\Mm) \leq l_{uR[F^{er}]}(\Mm) = l_{uR[F]}(\Mm) \leq l_{D_R}(\Mm),
\]
and furthermore, all these lengths are finite.
\end{proposition}
\begin{proof}
    By the above Theorem \ref{thm.Drfinite} the length of $\Mm$ as a $D_R$--module
    is finite. The inclusion of categories
    \[
        uR[F^e]\text{--mod} \subseteq uR[F]\text{--mod} \subseteq
        D_R\text{--mod}
    \]
    shows the claimed inequality of lengths. The equality follows,
    since a chain of finitely many $R[F]$--modules is a chain of
    $R[F^{er}]$--modules for some $r > 0$.
\end{proof}
More interesting is the case when one asks what happens after extending the
perfect field $k \subseteq R$ over which $R$ is essentially of finite type. If
$K$ is such an extension field of $k$, then we denote $R_K = K \tensor_k R$.
Clearly, tensoring with $K \tensor_k \usc$ is a functor from $uR[F^e]$--mod to
$uR_K[F^e]$--mod, from $uR[F]$--mod to $uR_K[F]$--mod and from $D_R$--mod to
$D_{R_K}$--mod, preserving finite generation. Furthermore, this functor can
only increase the length, since it is faithful.
\begin{definition}
    Let $\Mm$ be a finitely generated unit $R[F^e]$--module.
    The \emph{geometric length} of $\Mm$ is defined as the length of $K
    \tensor_k \Mm$, where $K$ is an algebraically closed field containing
    $k$. This notion of geometric length applies to all the lengths
    introduced above, and is denoted by $\overline{l}_*$, where $*$ is
    the appropriate category.
\end{definition}
For this to make sense, one has to show that the length does not depend on the
chosen algebraically closed field $K$. For the categories of finitely generated
unit $R[F^e]$--modules and $R[F]$--modules, this is a consequence of the next
proposition. The case of $D_R$--submodules of unit $R[F]$--modules is treated
in the next subsection as a byproduct of Theorem \ref{thm.LengthBigK}.
\begin{proposition}\label{prop.IndepField}
    Let $R$ be essentially of finite type over an algebraically closed
    field $k$ and let $K$ be an algebraically closed extension
    field. If $\Mm$ is a simple finitely generated unit $R[F^e]$--module
    (resp.\ $R[F]$--module) then $K \tensor_k \Mm$ is also simple as
    a unit $R_K[F^e]$--module (resp.\ $R_K[F]$--module).
\end{proposition}
\begin{proof}
The proof of this is a fairly standard argument using the Nullstellensatz and
generic flatness.\footnote{Thanks go to Brian Conrad for suggesting that these
techniques could be applied here.} In order to apply these techniques, one has
to work with finitely generated $R$--modules. Thus the trick consists of using
the roots of the unit $R[F^e]$--modules in question. These are finitely
generated $R$--modules.

Let $\beta: M \to \F[e]M$ be a root of $\Mm$. Then $\beta_K: M_K=K \tensor_k M
\to \F[e]_{R_K}M_K$ is a root of $\Mm_K=K \tensor_k \Mm$. Let $\Nn$ be a
nonzero unit $R[F^e]$--submodule of $\Mm_K$. Then the restriction of $\beta_K$
to $N_K = M_K \cap \Nn$ is a root of $\Nn$. Since both $M_K$ and $N_K$ are
finitely presented, one can find free presentations
\[
    R_K^{\oplus n'} \to R_K^{\oplus m'} \to M_K \to 0\quad\text{and}\quad   R_K^{\oplus n} \to R_K^{\oplus m} \to N_K \to
    0.
\]
The root morphism $\beta_K$ as well as the inclusion $i_K: N_K \to M_K$ can be
extended to a commutative diagram:
\begin{equation}\label{dia.K}
\begin{split}
\xymatrix{
    {R_K^{\oplus n'}} \ar[r] & {R_K^{\oplus m'}}\ar[r] & {M_K}\ar[r] & 0 \\
    {R_K^{\oplus n}} \ar[r]\ar[d]\ar[u] & {R_K^{\oplus m}}\ar[r]\ar[d]\ar[u] & {N_K}\ar[r]\ar_{i_K}@{^(->}[u]\ar^{\beta_K}@{^(->}[d] & 0  \\
    {\F[e]R_K^{\oplus n}} \ar[r]\ar@{=}[d] & {\F[e]R_K^{\oplus m}}\ar[r]\ar@{=}[d] &
    {\F[e]N_K}\ar[r]& 0  \\
    {R_K^{\oplus n}} \ar[r] & {R_K^{\oplus m}}
}
\end{split}
\end{equation}
The third row is just the Frobenius functor $\F[e]_{R_K}$ applied to the second
row. The equal signs indicate the natural unit $R[F^e]$--structure on $R$.
Thus, the six maps of free $R_K$--modules constituting the left big rectangle
completely determine the root morphisms $\beta_K$ and the inclusion $i_k:N_K
\subseteq M_K$. Since these are maps of finitely generated free $R_K$--modules,
they are defined over a finitely generated $k$--algebra $A \subseteq K$: just
adjoin to $k$ all the coefficients of the matrices representing these maps.
Then these maps are, in fact, maps of free $R_A = (A \tensor_k R)$--modules. We
get the corresponding diagram of $R_A$--modules \emph{defining} $M_A$ and $N_A$
as the cokernels of these maps of free modules (clearly, $M_A \cong A \tensor_k
M$):
\begin{equation}\label{dia.A}
\begin{split}
\xymatrix{
    {R_A^{\oplus n'}} \ar[r] & {R_A^{\oplus m'}}\ar[r] & {M_A}\ar[r] & 0 \\
    {R_A^{\oplus n}} \ar[r]\ar[d]\ar[u] & {R_A^{\oplus m}}\ar[r]\ar[d]\ar[u] & {N_A}\ar[r]\ar^{\beta_A}[d]\ar_{i_A}[u] & 0  \\
    {\F[e]R_A^{\oplus n}} \ar[r]\ar@{=}[d] & {\F[e]R_A^{\oplus m}}\ar[r]\ar@{=}[d] &
    {\F[e]N_A}\ar[r]& 0  \\
    {R_A^{\oplus n}} \ar[r] & {R_A^{\oplus m}}
}
\end{split}
\end{equation}
Since $\F[e]$ is right exact, the cokernel of the third row is in fact
$\F[e]_{R_A}N_A$ as indicated. The map $\beta_A$ and the map $i_A$ are the ones
induced on cokernels. By further enlarging $A$ one can assume that all the
modules involved (and especially the kernels of $\beta_A$ and $i_A$) are free
$A$--modules by generic flatness, \cite{Eisenbud.CommAlg}, Theorem 14.4.
Obviously, tensoring this diagram with $K$ over $A$, one gets back diagram
(\ref{dia.K}). Since we chose $A$ such that the kernel of $\beta_A$ is free, it
follows that this kernel must be zero, since it is zero after tensoring with
$K$. Similarly, the kernel of $i_A$ is also zero. Thus $\beta_A$ and $i_A$ are
injective. Thus we can think of $\beta$ as the root of a finitely generated
unit $R_A[F^e]$--submodule $\Nn_A$ of $\Mm_A=A \tensor_k \Mm$.

Now, let $m$ be a maximal ideal of $A$. Then, reducing mod $m$, we see that
$N_k
\defeq A/m \tensor_A N_A$ is a root of the unit $R[F^e]$--submodule $\Nn_k
\defeq A/m \tensor_A \Nn_A$ of $\Mm$ (here we used the Nullstellensatz and the
algebraic closedness of $k$ to conclude that $A/m \tensor_A R_A \cong R$ since
$A/m = k$). Again, by ensuring that $N_A$ is free over $A$ it follows that
$N_k$ and therefore $\Nn_k$ is nonzero. Since $\Mm$ is simple it follows that
$\Nn_k=\Mm$. But this implies that $N_k=M$, and therefore $N_A=M_A$. Thus
$N_K=M_K$ and consequntly $\Nn_K=\Mm_K$. This implies that $\Mm_K$ is a simple
unit $R_K[F^e]$--module.

The case of $\Mm$ being a simple $R[F]$--module follows easily. $\Mm$ is a
simple $R[F^e]$--module for infinitely many $e$. Thus $\Mm_K$ is a simple
$R_K[F^e]$--module for such $e$, and therefore also a simple $R_K[F]$--module.
\end{proof}
As a corollary of the proof of the proposition, one gets
\begin{corollary}
    Let $R$ be essentially of finite type over a perfect field $k$. Let $\Mm$
    be a finitely generated unit $R[F^e]$--module. Then there is a finite
    algebraic extension field $k'$ of $k$ such that $\bar{l}_{uR[F^e]}(\Mm)
    = l_{R_{k'}[F^e]}(\Mm_{k'})$, and similarly for the unit $R[F]$--module
    length.
\end{corollary}
\begin{proof}
After tensoring with the algebraic closure $K$, one has a finite sequence of
unit $R_K[F^e]$--submodules of $\Mm_K$:
\[
    0 \subseteq \Mm^1 \subseteq \Mm^2 \subseteq \ldots \subseteq \Mm^s =
    \Mm_K,
\]
where $s$ is the geometric $R[F^e]$--module length of $\Mm$. Similarly as in
the proof of the last proposition, \ie turning to presentations of the roots,
we see that all these $M^i$'s arise as $K \tensor_k M^i_{k'}$, where $M^i_{k'}$
are unit $R_K[F^e]$--submodules of $M_{k'}$ ($k'$ takes the role of $A$ in the
last proof). This shows that the unit $R_{k'}[F^e]$--module length of $M_{k'}$
is equal to the geometric unit $R[F^e]$--module length of $\Mm$.

The case of unit $R[F]$--module length is achieved by looking at
$R[F^e]$--module lengths for sufficiently big $e$.
\end{proof}

\subsection{Geometric $F$--length is geometric
$D_R$--length}\label{sec.RF=DR} The fact that the geometric $D_R$--module
length is well defined will be proved by showing that for sufficiently huge
fields $K$, the $D_{R_K}$--module length is equal to the $R_K[F]$--module
length. Then, if $R$ is essentially of finite type over the algebraically
closed field $k$, we have the chain of inequalities
\[
    \bar{l}_{uR[F]} \leq l_{D_R} \leq l_{D_{R_K}} = l_{uR_K[F]}.
\]
Therefore equality prevails everywhere, and the geometric $D_R$--module length
is also well defined. The crucial step is therefore the following theorem.
\begin{theorem}\label{thm.LengthBigK}
    Let $k$ be algebraically closed and of strictly bigger cardinality than
    a $k$--basis of $R$. Then a simple finitely generated unit $R[F]$--module $\Mm$ is simple as a
    $D_R$--module.
\end{theorem}
\begin{proof}
By Corollary \ref{cor.IsotypicAreSub}, $(\Mm, \theta^e_{\Mm}, F^e_{\Mm})$ is
$D_R$--isotypic. Hence, as a $D_R$--module, $\Mm$ is isomorphic to $V \tensor_k
\Nn$ for a simple $D_R$--module $\Nn$ and a finite-dimensional
$k$--vector-space $V$. The cardinality assumption together with the algebraic
closedness of $k$ implies (Lemma \ref{lem.EndNisAlgebraic}) that
$\End_{D_R}(\Nn)=k$. Thus we are in the situation of Proposition
\ref{prop.MreducetoV}, and therefore obtain a unit $k[F^e]$--structure $F^e_V$
on $V$ such that the unit $R[F]$--submodules of $\Mm$ are in one-to-one
correspondence with the $k[F]$--submodules of $V$. In particular, $\Mm$ is
simple as a unit $R[F]$--module if and only if $V$ has no nontrivial
$F^e_V$--stable submodules. The following Proposition \ref{prop.dioudonne} of
Dieudonn{\'e} then shows that $V$ is one-dimensional. Therefore $\Mm \cong \Nn$ and
$\Mm$ is $D_R$--simple.
\end{proof}
The last technical aid is Dieudonn{\'e} \cite{DieudonLie}, Proposition 3, page 233;
the statement of which we recall.
\begin{proposition}\label{prop.dioudonne}
   Let\/ $k$ be an algebraically closed field and\/ $V$ a finite-dimensional\/
   $k[F^e]$--vector-space. If the Frobenius\/ $F^e$ acts injectively on\/ $V$, then\/
   $V$ has a basis consisting of\/ $F^e$--fixed elements of\/ $V$.
\end{proposition}
As a corollary of the last theorem we get the equality of geometric
$R[F]$--module length and $D_R$--module length.
\begin{corollary}
    Let $R$ be regular, essentially of finite type over a field $k$,
    and let $\Mm$ be a finitely generated unit $R[F]$--module. Then
    the geometric $D_R$--module length of $\Mm$ is well defined  and
    equal to the geometric unit $R[F]$--module length of $\Mm$.
\end{corollary}
\begin{proof}
Since the geometric unit $R[F]$--module length is well defined (Proposition
\ref{prop.IndepField}),
\[
\bar{l}_{uR[F]}(\Mm)=l_{uR_K[F]}(\Mm_K)
\]
with $K$ any algebraically closed extension field of $k$. Given a maximal
filtration of $\Mm_K$ as a unit $R_K[F]$--module,
\[
    0 \subseteq \Mm^1 \subseteq \ldots \subseteq \Mm^l = \Mm_K,
\]
with simple unit $R_K[F]$--module quotients, if $K$ is sufficiently big
(uncountable), the last theorem shows that all the quotients are
$D_{R_K}$--simple. Therefore the geometric $D_R$--module length of $\Mm$ is
equal to $l=\bar{l}_{uR[F]}$, and thus it too is well defined (independent of
$K$).
\end{proof}

\section{Examples}
An example of a simple $R[F]$--module that is not simple as a
$D_R$--module is given. This example is constructed as a free
$R$--module $\Mm$ of rank $2$. Let the action of the Frobenius
$F^e$ with respect to some basis $\Aa=(e_1,e_2)$ be represented by
the matrix
\[
    A= \begin{pmatrix} 0 & 1 \\ 1 & x \end{pmatrix}
\]
for some element $x \in R$, \ie for $v = v_1e_1+v_2e_2$ the action of $F^e$ is
given by
\[
    v=\begin{pmatrix} v_1 \\ v_2 \end{pmatrix} \mapsto A\begin{pmatrix} v^q_1 \\ v^q_2
    \end{pmatrix}.
\]
The choice of basis $\Aa$ also induces a natural $D_R$--structure
on $\Mm \cong R\oplus R$ by acting componentwise on the direct
summands. Unless otherwise specified, this is the $D_R$--module
structure on $\Mm$ we have in mind.

The matrix $A_r$ representing the $r$th power of this Frobenius
action $F^e$ with respect to this basis is given by
\[
    A_r= AA^{[q]}\cdots A^{[q^{r-1}]},
\]
where the square brackets $[q]$ raise each coefficient of the matrix to its
$q$th power. Equivalently, $A_r$ can be described inductively by the equation
$A_r = A_{r-1}A^{[q^{r-1}]}$, which translates into an inductive formula for
the coefficients of $A_r$. One has
\begin{equation}
    A_r=\begin{pmatrix} a^q_{r-2} & a^q_{r-1} \\
                       a_{r-1}   & a_r
        \end{pmatrix},
\end{equation}
where $a_r=a_{r-2}+a_{r-1}x^{q^{r-1}}$ with $a_{-1}=0$ and
$a_0=1$. So, for example, this formula computes
\[
    A_1 = \begin{pmatrix} {0} &{1} \\ {1} &{x} \end{pmatrix}\ \text{ and }\
    A_2 = \begin{pmatrix} {1} &{x^q} \\ {x} &{x^{q+1}+1}
\end{pmatrix},
\]
which can be easily verified by hand. We prove these assertions by induction:
\begin{equation*}\begin{split}
    A_r = A_{r-1}A^{[q^{r-1}]} &=
      \begin{pmatrix}
            a^q_{r-2} & a^q_{i-3}+a^q_{r-2}x^{q^{r-1}}  \\
            a_{i-1}   & a_{a-2}+a_{i-1}x^{q^{r-1}}
      \end{pmatrix} \\  &= \begin{pmatrix}
                            a^q_{r-2} & (a_{r-3}+a_{r-2}x^{q^{r-2}})^q \\
                            a_{r-1}   & a_r
                     \end{pmatrix} = \begin{pmatrix} a^q_{2-1} & a^q_{r-1} \\
                                                    a_{r-1}   & a_r
                                    \end{pmatrix},
\end{split}\end{equation*}
where the base case is by the initial condition of the recursion for $a_r$.
Furthermore, we note that, thinking of $a_r$ as polynomials in $x$, the degree
in $x$ of $a_r$ is $\deg a_r = 1 + q + q^2 +\ldots+q^{r-1}$. Again an induction
argument shows this nicely:
\begin{equation*}\begin{split}
    \deg a_r &= \max\{\deg a_{r-2}, \deg a_{r-1}+q^{r-1}\} \\
             &= \max\{1+q+\ldots q^{r-3},1+q+\ldots+q^{r-2}+q^{r-1}\} \\
             &= 1+q+\ldots+q^{r-2}+q^{r-1},
\end{split}\end{equation*}
and for the start of the induction we just recall that $a_{-1}=0$ and $a_0=1$.
This setup is used as the basis for the following examples.

\subsubsection{A simple \protect{$R[F^1]$--module that is not
$D_R$}--simple} Let $R=\FF_3$ and $x=1$. Since the Frobenius $F$
is the identity on ${\FF_3}$, Frobenius actions are just linear
maps. The linear map represented by
\[
   A = \begin{pmatrix} 0 & 1 \\ 1 & 1 \end{pmatrix}
\]
is not diagonalizable, since its characteristic polynomial $P_A(t)=t(t-1)+1$ is
irreducible over $\FF_3$. Thus $\Mm$ is a simple ${\FF_3}[F^1]$--module. Since
${\FF_3}$ is perfect, $D_{\FF_3}={\FF_3}$, and thus $\Mm$ is not simple as a
$D_{\FF_3}$--module, since it is a free ${\FF_3}$--module of rank 2. Note that,
since $\End_{\FF_3}(\Mm)$ is finite, some power of $F$ will be the identity on
$\Mm$. An easy calculation shows that $F^4=-\id_M$, and therefore $\Mm$ is not
simple as a ${\FF_3}[F^4]$--module.

On the other hand, after adjoining a root $\alpha$ of the polynomial
$P_A=t(t-1)+1$ to $\FF_3$ we get the extension field $k=\FF_3(\alpha)$.
Clearly, $k\tensor_{\FF_3} \Mm$ is no longer simple as a $k[F^1]$--module.

\subsubsection{\protect{A $D_R$--submodule that is not an
$R[F]$}--submodule}\label{ex.DRnorRFinfty} Now let $R$ be a ring containing an
infinite perfect field $k$ with $x \in k$ transcendental over the prime field
$\FF_p$. With this $x$ the matrix $A$, in fact, represents a $D_R$--linear map
of $\Mm$ (this is because differential operators $D_R$ are linear over any
perfect subring of $R$). Since $a_r$ is a nonzero polynomial in $x$ with
coefficients in $\FF_p$, and since $x$ is transcendental over the prime field,
$a_r$ is a nonzero element of $k$. This implies that, for example, $Re_1$ is
not stable under any power of $F^e$, since this would be equivalent to the
matrix $A_r$ having a zero entry in the bottom left corner. But this entry is
$a_{r-1}$, which, as we just argued, is nonzero. Thus $Re_1$ is a (simple)
$D_R$--submodule that is not an $R[F]$--submodule of $M$.

\subsection{A simple unit \protect{$R[F]$}--module that is not
$D_R$--simple} Now we come to the main example of a simple $R[F]$--module that
is not $D_R$-simple. This will show that the unit $R[F]$--module and
$D_R$--module lengths of a unit $R[F]$--module can be different. It provides a
counterexample to Lyubeznik's Remark 5.6a in \cite{Lyub}, where he speculates
that simple $D_R$--submodules of simple unit $R[F]$--modules are in fact
$R[F]$--submodules.

With $A$ and $\Mm$ as before, let $R=k(x)^{1/p^\infty}$ be the perfect closure
of $k(x)$, where $x$ is a new variable and $k$ is perfect. As before, we define
the Frobenius action $F^e$ on $\Mm$ as being represented by $A$, \ie given by
application of
\[
    F^e(\usc) = \begin{pmatrix} 0 & 1 \\ 1 & x \end{pmatrix}(\usc)^{[p^e]}.
\]
with respect to the basis $\Aa$. Since $R$ is perfect $R=D_R$, and consequently
this (and every) unit $R[F]$--structure is compatible with the
$D_R$--structure.

To show that $(\Mm,F^e)$ is simple as a unit $R[F]$--module, we show that $M$
is a simple unit $R[F^{er}]$--module for all $r$. We proceed in 3 steps:
\begin{descrip}{Step 1\phantom{a}}
\item[\textsc{Step 1}]
    With respect to the basis $\Aa$ the action $F^{er}$ is represented by
    the matrix $A_r$. We change the basis appropriately to $\Bb=(f_1,f_2)$
    such that the representing matrix $B_r$ of $F^{er}$ with respect to
    the basis $\Bb$ is ``nice''; by this we mean that
    \[
     B_r = \begin{pmatrix} 0 & s_r \\ 1 & t_r \end{pmatrix}
    \]
    for some $s_r, t_r \in \FF_p[x]$.
\item[\textsc{Step 2}]
    Assuming that there is a $v \in M$ such that
    $F^{er}(v) = \lambda v$ yields a monic algebraic equation which gives an
    algebraic equation for $x$.
\item[\textsc{Step 3}]
    By transcendence of $x$ this equation must be zero, and we discriminate
    two cases to arrive at a contradiction. One is
    treated by a degree argument, the other by differentiation.
\end{descrip}

Let us begin with \textsc{Step 1}: The basis that will lead to the
matrix $B_r$ of the desired shape is $f_1 = e_1$ and $f_2 =
a^q_{r-2}e_1 + a_{r-1}e_2$. Thus the matrix responsible for the
base change from $\Aa$ to $\Bb$ is
\[
    C_r \defeq \begin{pmatrix} 1 & a^q_{r-2} \\ 0 & a_{r-1} \end{pmatrix}.
\]
With respect to the new basis $\Bb$ the Frobenius action $F^{er}$ is
represented by the matrix $B_r = C_r^{-1}A_r C_r^{[q^{r}]}$. This can be
checked by hand; a more thorough discussion of Frobenius actions on free
$R$--modules under change of basis can be found in \cite{Manuel.PhD}. To
determine $s_r$ and $t_r$ we explicitly calculate $B_r$:
\begin{equation*}
\begin{split}
B_r &= C_r^{-1}A_r C_r^{[q^r]} = C_r^{-1}
\begin{pmatrix}
    {a_{r-2}^q}&{a_{r-1}^q}\\
    {a_{r-1}}&{a_r}
\end{pmatrix}
\begin{pmatrix}
    {1}&{a_{r-2}^{q^r+q}}\\
    {0}&{a_{r-1}^{q^r}}
\end{pmatrix} \\
&=\frac{1}{a_{r-1}}
\begin{pmatrix}
    {a_{r-1}}&{-a_{r-2}^q} \\
    {0}& {1}
\end{pmatrix}
\begin{pmatrix}
    {{a_{r-2}^q}} & {{a_{r-2}^q a_{r-2}^{q^r+q} + a_{r-1}^q a_{r-1}^{q^r}}} \\
    {{a_{r-1}}} & {{a_{r-1}a_{r-2}^{q^r+q} + a_r a_{r-1}^{q^r}}}
\end{pmatrix} \\
%&= \frac{1}{a_{r-1}}   % This might not be necessary...
%\begin{pmatrix}
%    {a_{r-1}a_{r-2}^q-a_{r-2}^q a_{r-1}}&
%    \scriptstyle{{a_{r-1}a_{r-2}^{q^r+2q}+a_{r-1}^{q^r+q+1}-a_{r-1}a_{r-2}^{q^r+2q}-a_ia_{r-1}^{q^r}a_{r-2}^q}} \\
%%\end{pmatrix}\\
&=
\begin{pmatrix}
    {0} &{-a_{r-1}^{q^r-1}\det A_r}\\
    {1} &{a_{r-2}^{q^r+q} + a_ra_{r-1}^{q^r-1}}
\end{pmatrix} \\
&=
\begin{pmatrix}
    {0}&{(-1)^{r-1}a_{r-1}^{q^r-1}} \\
    {1}&{a_{r-2}^{q^r+q} + a_ra_{r-1}^{q^r-1}}
\end{pmatrix}.
\end{split}
\end{equation*}
Besides index juggling skills, one only needs the equation $\det A_r = (-1)^r$,
which follows from the  recursive definition of $A_r$ and the fact that $\det
A=-1$. We can read off the desired expressions for $s_r$ and $t_r$:
\[
    s_r=(-1)^{r-1}a_{r-1}^{q^r-1} \text{\quad and \quad}
t_r=a_{r-2}^{q^r+q} + a_r a_{r-1}^{q^r-1}.
\]
Note that both are in $\FF_p[x]$, since $a_i \in \FF_p[x]$. We now work over
this new basis $\Bb$ and start with \textsc{Step 2}. Assume we have
$v=(\alpha,1)^t$ such that $F^r(v)=\lambda v$. One reduces from a general
$v=(\alpha,\beta)^t$ to this case by dividing by $\beta$; \ref{ex.DRnorRFinfty}
ensures that $\beta$ is not zero. Considering $\lambda v = B_rv^{[q^r]}$, this
yields two equations:
\begin{equation*}
\begin{split}
    s_r &=  \lambda \alpha, \\
    \alpha^{q^r}+t_r &= \lambda.
\end{split}
\end{equation*}
Substituting the latter in the former, we get a monic algebraic equation for
$\alpha$ with coefficients in $k[x]$:
\begin{equation*}
    \alpha^{q^r+1} + t_r \alpha - s_r = 0
\end{equation*}
Since $k[x]^{1/p^{\infty}}$ is integrally closed in $R$, we conclude that
$\alpha \in k[x]^{1/p^{\infty}}$. Choose $t$ minimal such that $\alpha \in
k[x^{1/p^t}]$. Then $\beta(x)=\alpha^{p^t}$ is in $k[x]$, and not a $p$th power
unless $t=0$. Taking the $p^t$th power of the last equation, we get
\begin{equation}\label{eqn.betar}
    \beta^{q^r+1} + t_r^{p^t} \beta - s_r^{p^t} = 0,
\end{equation}
which is an algebraic relation for $x$ with coefficients in $k$. Thus it is
constant zero by transcendence of $x$. For \textsc{Step 3} we distinguish the
following two cases:
\begin{descrip}{$t > 0\ :$}
\item[$t > 0$ :] Differentiating \eqnref{eqn.betar} with respect to $x$, we
get:
\[
    \beta^{q^r}\diff \beta + t_r^{p^t} \diff \beta = 0.
\]
Since we chose $\beta$ not to be a $p$th power, its derivative is nonzero. Thus
we can divide the above by $\diff \beta$ and get $\beta^{q^r}=-t_r^{p^t}$.
Substituting this back into \eqnref{eqn.betar}, we get $s_r^{p^t}$=0. But this
is a contradiction, since $s_r=\pm a_{r-1}^{q^r-1} \neq 0$.

\item[$t=0$ :] For this we have to determine the degrees of the terms in
equation \eqnref{eqn.betar}. As observed earlier,
$\deg(a_r)=1+q+\ldots+q^{r-1}$. Therefore,
\begin{equation*}
  \begin{split}
  \hbox{\phantom{$t=0$}}\deg{s_r}&=(q^r-1)\deg(a_{r-1}) \\
                                 &=-1-q-\ldots-q^{r-2}+q^r+\ldots+q^{2r-2}, \\
                        \deg{t_r}&\leq \max\{\, \deg(a_{r-2}^{q^r+q}),\ \deg({a_ia_{r-1}^{q^r-1}})\, \} \\
                                 &=q^{r-1}+\ldots+q^{2r-2}.
  \end{split}
\end{equation*}
In fact, equality prevails in the last inequality, since the two entries in the
$\max$ are different (the second is always bigger). To be precise,
\begin{equation*}
\begin{split}
\deg(a_{r-2}^{q^r+q}) &=  (q^r+q)(1+\ldots+q^{r-3}) \\
                      &= q+\ldots+q^{r-2}+q^r+\ldots+q^{2r-3}, \\
\hbox{\phantom{$t=0$}}\deg({a_ia_{r-1}^{q^r-1}}) &= 1+q+\ldots+q^{r-1}-1-\ldots \\
                                                 &\qquad \qquad \qquad -q^{r-2}+q^r+\ldots+q^{2r-2} \\
                           &= q^{r-1}+\ldots+q^{2r-2}
\end{split}
\end{equation*}
Since $q^{r-1} > 1+q+\ldots+q^{r-2}$, we see that the second line is in fact
strictly bigger than the first. Thus the degree of $s_r$ is strictly smaller
than the degree of $t_r$, and therefore the first two terms of
\eqnref{eqn.betar} must have the same degree. If we denote the degree of
$\beta$ by $n$, we get
\[
    (q^r+1)n = \deg(t_r)+n,
\]
and after dividing by $q^{r-1}$ this simplifies to
\[
    1 = qn-q-q^2-\ldots-q^{r-1}.
\]
The right side is divisible by $q$, but the left side certainly is not. This is
a contradiction.
\end{descrip}
This finishes the proof that $\Mm$ is a simple $R[F^{er}]$--module for all
$r>0$. Thus $\Mm$ is a simple $R[F]$--module, but $M$ is not simple as a
$D_R$--module since every one dimensional $R$--subspace is a nontrivial
$D_R$--submodule.

Let $R'=R(\alpha)$ be the field one obtains by adjoining to $R$ a
root $\alpha$ of the polynomial $P(t)=t^{p^2}+xt^p-t$. Then
\[
    F(\begin{pmatrix}\alpha^p \\ \alpha\end{pmatrix}) =
     \begin{pmatrix} 0 & 1 \\ 1 & x\end{pmatrix}
     \begin{pmatrix}\alpha^{p^2} \\ \alpha^p\end{pmatrix} =
     \begin{pmatrix}\alpha^p \\ \alpha^{p^2}+x\alpha^p\end{pmatrix}=
     \begin{pmatrix}\alpha^p \\ \alpha\end{pmatrix},
\]
which shows that the element $(\alpha^p, \alpha)^t$ of $\Mm' = R'
\tensor_R \Mm$ is fixed by $F$. Thus $\Mm'$ is not simple as an
$R'[F]$--module.

\subsubsection{Examples over the polynomial ring} So far the
examples were over a field. Starting with these examples it is not hard to
obtain equivalent examples over higher-dimensional rings. For this let
$(V,F^e)$ be the simple unit $K[F]$--module of the last example ($K =
\FF_p(x)^{1/p^\infty}$). Let $R$ be a regular $K$--algebra, essentially of
finite type over $K$ (\eg $R=K \tensor_{\FF_p} \FF_p[\xn]=K[\xn]$). Let $\Mm =
R \tensor_K V$, where we point out that $R$ is a simple $D_R$-module and
$\End_{D_R}(R)=K$. Then $\Mm$ carries a natural $R[F^e]$--structure defined by
$F^e_\Mm(n \tensor v) = n^{p^e} \tensor F^e(v)$. This is exactly the situation
of Proposition \ref{prop.MreducetoV}, and it follows that the unit
$K[F^e]$--submodules of $V$ are in one-to-one correspondence with the unit
$R[F^e]$--submodules of $\Mm$. Thus, the simplicity of $V$ as a unit
$K[F]$--module implies that $\Mm$ is a simple unit $R[F]$--module. Clearly,
since $V$ is not $D_K$--simple, $\Mm$ is not $D_R$--simple.

\end{document}